\documentclass[12pt,a4paper,reqno]{amsart}
\usepackage{graphics,epic}
\usepackage{amsmath,amssymb, amsthm}
\usepackage[all,2cell]{xy}

\newtheorem{theorem}{Theorem}[section]
\newtheorem*{theorem*}{Theorem}

\newtheorem{lemma}[theorem]{Lemma}
\newtheorem{proposition}[theorem]{Proposition}
\newtheorem{corollary}[theorem]{Corollary}

\newtheorem*{conjecture*}{Conjecture}

\newtheorem{remark}[theorem]{Remark}

\renewcommand{\hat}[1]{\widehat{#1}}

\newcommand{\f}{\mathfrak{f}}

\newcommand{\id}{{\rm id}}

\newcommand{\End}{{\rm End}\,}
\newcommand{\Res}{{\rm Res}\,}

\newcommand{\Z}{\mathbb{Z}}

\newcommand{\C}{\mathbb{C}}

\def\Res{{\rm Res}}
\def\wt{{\rm wt}}
\def\C{{\mathbb C}}

\def\Z{{\mathbb Z}}

\def\1{{\bf 1}}
\def\tr{{\rm tr}}
\def \End{{\rm End}}

\def \Ind{{\rm Ind}}
\def \pf{\noindent {\bf Proof: \,}}
\def\theequation{5.\arabic{equation}}

\def \h{\mathfrak{h}}

\def \w{\omega}
\def \g{\mathfrak{g}}

\begin{document}

\title[Module category and $C_2$-cofiniteness of affine VOSAs]{Module category and $C_2$-cofiniteness of affine vertex operator superalgebras}
\author{Chunrui Ai}
\address{Chunrui Ai, School of Mathematics and Statistics, Zhengzhou University, Henan 450001, China}
\email{aicr@zzu.edu.cn}

\thanks{C. Ai was supported by China NSF grant 11701520;  X. Lin was supported by China NSF grant
11801419 and the starting research fund from Wuhan University (No. 413000076)}
\author{Xingjun Lin}
\address{Xingjun Lin,  School of Mathematics and Statistics, Wuhan University, Wuhan 430072, China.}
\email{linxingjun88@126.com}

\begin{abstract}
 In this paper, we investigate the Lie algebra structures of weight one subspaces of $C_2$-cofinite vertex operator superalgebras. We also show that for any positive integer $k$, vertex operator superalgebras $L_{sl(1|n+1)}(k,0)$ and $L_{osp(2|2n)}(k,0)$ have infinitely many irreducible admissible modules. As a consequence, we give a proof of the fact that $L_{\g}(k,0)$ is $C_2$-cofinite if and only if $\g$ is either a simple Lie algebra, or $\g=osp(1|2n)$, and $k$ is a nonnegative integer. As an application, we show that $L_{G(3)}(1,0)$ is a vertex operator superalgebra such that the category of $L_{G(3)}(1,0)$-modules is semisimple but $L_{G(3)}(1,0)$ is not $C_2$-cofinite.
\end{abstract}
\maketitle
\section{Introduction }
\def\theequation{1.\arabic{equation}}
\setcounter{equation}{0}

Let $\g$ be a finite dimensional simple Lie superalgebra with a nondegenerate even supersymmetric invariant bilinear form, $L_{\g}(k,0)$ be the simple vertex operator superalgebra associated with $\g$ \cite{K2}.  In case that $\g$ is a finite dimensional simple Lie algebra and  $k$ is a positive integer, it was proved that $L_{\g}(k,0)$ is $C_2$-cofinite \cite{DLM2}, \cite{Z}. Moreover, in case that $\g$ is a finite dimensional simple Lie algebra, it was proved in \cite{DM} that $k$ must be a positive integer if $L_{\g}(k,0)$ is $C_2$-cofinite. It has also been known for many years that $L_{\g}(k,0)$ is $C_2$-cofinite if and only if $\g$ is either a simple Lie algebra, or $\g=osp(1|2n)$, and the $\hat{ \g}$-module $L_{\g}(k,0)$ is integrable (see Section 0.4 of  \cite{GK}). One of our motivations is to give a proof of this result.

Our first main result is about the Lie algebra structures of weight one subspaces of $C_2$-cofinite vertex operator superalgebras. Explicitly, let  $V$ be a  vertex operator superalgebra which is of strong CFT type and $C_2$-cofinite. It is known \cite{B} that  $V_1$ has a Lie superalgebra structure. Assume that $\g$ is a subalgebra of $V_1$ and that $\g$ is a basic simple Lie superalgebra, then we show that $\g$ must be isomorphic to $osp(1|2n)$, $sl(1|n+1)$ or $osp(2|2n), n\geq 1$ (see Theorem \ref{inte}). Our second main result is about admissible modules of vertex operator superalgebras $L_{sl(1|n+1)}(k,0)$ and $L_{osp(2|2n)}(k,0)$. We show that $L_{sl(1|n+1)}(k,0)$ and $L_{osp(2|2n)}(k,0)$ have infinitely many irreducible admissible modules (see Theorems \ref{non1}, \ref{non2}). As an application of these results, we give a proof of the fact that $L_{\g}(k,0)$ is $C_2$-cofinite if and only if $\g$ is either a simple Lie algebra, or $\g=osp(1|2n)$, and $k$ is a nonnegative integer (see Theorem \ref{main}).

Another motivation of this work comes from the study of the relationship between rationality and $C_2$-cofiniteness of vertex operator algebras. For a vertex operator algebra $V$, we could consider categories of weak $V$-modules, admissible $V$-modules and ordinary $V$-modules \cite{DLM3}. A vertex operator algebra $V$ is called regular if the category of weak $V$-modules is semisimple, and $V$ is called rational if the category of admissible $V$-modules is semisimple. In \cite{Z}, Zhu conjectured that rational vertex operator algebras are $C_2$-cofinite. It has been proved that regular vertex operator algebras are $C_2$-cofinite \cite{Li0}. However, Zhu's conjecture has not been proved. This motivates us to ask whether a vertex operator superalgebra $V$ is $C_2$-cofinite if the category of ordinary $V$-modules is semisimple.  Due to Kac and Wakimoto, it is known that the category of ordinary $L_{G(3)}(1,0)$-modules is semisimple (see Remark 8.2 of \cite{KWaki}). In this paper, we give a proof of this fact (see Theorem \ref{semisimple}). As a consequence, we show that $L_{G(3)}(1,0)$  is a vertex operator superalgebra such that the category of $L_{G(3)}(1,0)$-modules is semisimple but $L_{G(3)}(1,0)$ is not $C_2$-cofinite.

 Our work was also motivated by the study of coset vertex operator subalgebras of affine vertex operator algebras \cite{ACL}, \cite{FZ}. It is known \cite{W1}, \cite{W2} that coset vertex operator subalgebras of affine vertex operator algebras are closely related to $W$-superalgebras, which are vertex operator superalgebras obtained from affine vertex operator superalgebras by  quantum Drinfeld-Sokolov reduction \cite{KRW}.  We expect that our results are useful for studying $W$-superalgebras and coset vertex operator subalgebras of affine vertex operator algebras.

 The paper is organized as follows: In Section 2, we recall basic definitions about vertex operator superalgebras and basic facts about affine vertex operator superalgebras. In Section 3, we investigate the Lie algebra structures of weight one subspaces of $C_2$-cofinite vertex operator superalgebras. In Section 4, we show that $L_{sl(1|n+1)}(k,0)$ and $L_{osp(2|2n)}(k,0)$ have infinitely many irreducible admissible modules. As a consequence, it is proved that $L_{sl(1|n+1)}(k,0)$ and $L_{osp(2|2n)}(k,0)$ are not $C_2$-cofinite. In Section 5, we prove that $L_{osp(1| 2n)}(k, 0)$ is a $C_2$-cofinite vertex operator superalgebra. As a consequence, it is proved that $L_{\g}(k,0)$ is $C_2$-cofinite if and only if $\g$ is either a simple Lie algebra, or $\g=osp(1|2n)$, and $k$ is a nonnegative integer. In Section 6, we show that the category of ordinary $L_{G(3)}(1,0)$-modules is semisimple.

\section{Preliminaries}
\def\theequation{2.\arabic{equation}}
\setcounter{equation}{0}
\subsection{Basics} In this subsection, we recall from  \cite{FHL}, \cite{KWa}, \cite{Li3} basic notions about vertex operator superalgebras. Let $V=V_{\bar 0}\oplus V_{\bar 1}$ be any ${\Z}_2$-graded vector space, the element in $V_{\bar 0} $ (resp. $V_{\bar 1}$) is called  {\em even} (resp. {\em odd}).
We then define $[v]=i$ for any $v\in V_{\bar i}$ with  $i=0,1$.  A {\em vertex superalgebra} is a quadruple $(V,Y(\cdot, z),\1,D),$ where $V=V_{\bar 0}\oplus V_{\bar 1}$ is a ${\Z}_2$-graded vector space, ${\bf 1}$ is the {\em vacuum vector}
of $V$, $D$ is an endomorphism of $V$,   and $Y(\cdot, z)$ is a linear map
\begin{align*}
 Y(\cdot, z): V &\to (\End\,V)[[z,z^{-1}]] ,\\
 v&\mapsto Y(v,z)=\sum_{n\in{\Z}}v_nz^{-n-1}\ \ \ \  (v_n\in
\End\,V)
\end{align*}
satisfying the following axioms:  \\
(i) For any $u,v\in V,$ $u_nv=0$ for sufficiently large $n$;\\
(ii) $Y({\bf 1},z)=\id_{V}$;\\
(iii) $Y(v,z){\bf 1}=v+\sum_{n\geq 2}v_{-n}{\bf 1}z^{n-1},$ for any $v\in V$;\\
(iv) $[D, Y(v, z)]=Y(D(v), z)=\frac{d}{dz}Y(v, z)$;\\
(v) The {\em Jacobi identity} for ${\Z}_2$-homogeneous $u,v\in V$ holds,
\begin{align*}
\begin{array}{c}
\displaystyle{z^{-1}_0\delta\left(\frac{z_1-z_2}{z_0}\right)
Y(u,z_1)Y(v,z_2)-(-1)^{[u][v]} z^{-1}_0\delta\left(\frac{z_2-z_1}{-z_0}\right)
Y(v,z_2)Y(u,z_1)}\\
\displaystyle{=z_2^{-1}\delta
\left(\frac{z_1-z_0}{z_2}\right)
Y(Y(u,z_0)v,z_2)}.
\end{array}
\end{align*}
This completes the definition of a vertex superalgebra and we will denote the vertex superalgebra briefly by $V$. A vertex  superalgebra $V$ is called {\em $C_2$-cofinite} if $\dim V/C_2(V)<\infty$, where $C_2(V)=\langle u_{-2}v|u, v\in V \rangle$.

A vertex superalgebra $V$ is called a {\em vertex operator superalgebra} if there is  a distinguished vector $\omega$, which is called the {\em conformal vector} of $V$, such that the following two conditions hold: \\
 (vi) The component operators of  $Y(\w,z)=\sum_{n\in\Z}L(n)z^{-n-2}$ satisfy the Virasoro algebra
relation with {\em central charge} $c\in \C:$
\begin{align*}
[L(m),L(n)]=(m-n)L(m+n)+\frac{1}{12}(m^3-m)\delta_{m+n,0}c,
\end{align*}
and
$$L(-1)=D;$$
(vii) $V$ is $\frac{1}{2}\Z$-graded such that $V=\oplus_{n\in \frac{1}{2}\Z}V_n$,  $L(0)|_{V_n}=n$, $\dim(V_n)<\infty$ and $V_n=0$ for sufficiently small $n$.

\begin{remark}
We adopt the definition of a vertex operator superalgebra given in \cite{Li3}. The definition of a vertex operator superalgebra is different from that in \cite{DZ}, \cite{Li}, \cite{X}, where the additional condition $V_{\bar 0}=\oplus_{n\in \Z}V_n$, $V_{\bar 1}=\oplus_{n\in \frac{1}{2}+\Z}V_n$ is added.  However, we will see that affine vertex operator superalgebras do not satisfy the  condition $V_{\bar 0}=\oplus_{n\in \Z}V_n$, $V_{\bar 1}=\oplus_{n\in \frac{1}{2}+\Z}V_n$ (see Theorem \ref{affine}).
\end{remark}

A vertex operator superalgebra $V$ is called of {\em CFT} type if $V$ has the decomposition $V=\bigoplus_{n\geq 0}V_n$ with respect to $L(0)$ such that $\dim V_0=1$. We say that $V$ is of {\em strong CFT type} if $V$ satisfies the further condition $L(1)V_1=0$. If $v\in V_n$, the {\em conformal weight} $\wt v$ of $v$ is defined to be $n$.

\subsection{Zhu's algebra of vertex operator superalgebras} In this subsection, we recall from \cite{DZ1} some facts about the Zhu's algebra of vertex operator superalgebras. First,  a {\em $\Z$-graded vertex operator superalgebra} $V$ is a vertex operator superalgebra such that $V_n=0$ for $n\in \frac{1}{2}+\Z$.
For a $\Z$-graded vertex operator superalgebra $V$, a {\em weak  $V$-module} $M$ is a vector space equipped
with a linear map
\begin{align*}
Y_{M}(\cdot, z):V&\to (\End M)[[z, z^{-1}]],\\
v&\mapsto Y_{M}(v,z)=\sum_{n\in\Z}v_nz^{-n-1},\,v_n\in \End M
\end{align*}
satisfying the following conditions: For any ${\Z}_2$-homogeneous  $u\in V,\ v\in V,\ w\in M$ and $n\in \Z$,
$$u_nw=0 \text{ for sufficiently large } n;$$
$$Y_M(\1, z)=\id_M;$$
\begin{align*}
\begin{split}
&z_{0}^{-1}\delta\left(\frac{z_{1}-z_{2}}{z_{0}}\right)Y_{M}(u,z_{1})Y_M(v,z_{2})-(-1)^{[u][v]}z_{0}^{-1}\delta\left(
\frac{z_{2}-z_{1}}{-z_{0}}\right)Y_M(v,z_{2})Y_M(u,z_{1})\\
&\quad=z_{2}^{-1}\delta\left(\frac{z_{1}-z_{0}}{z_{2}}\right)Y_M(Y(u,z_{0})v,z_{2}).
\end{split}
\end{align*}

A weak $V$-module $M$ is called {\em admissible} if it is $\Z_{\geq 0}$-graded $$M=\oplus_{n\in \Z_{\geq 0}}M(n)$$ such that for homogeneous $v\in V$,
$$v_mM(n)\subseteq M(n+\wt v-m-1).$$
In this paper, we will assume that $M(0)\neq 0$ if $M\neq 0$. A $\Z$-graded vertex operator superalgebra $V$ is called {\em rational} if any admissible $V$-module is completely reducible.

For a $\Z$-graded vertex operator superalgebra $V$, an {\em (ordinary)  $V$-module} is a weak $V$-module $M$ which carries a $\C$-grading induced by the spectrum of $L(0)$, that is,  $M=\bigoplus_{\lambda\in\C}
M_{\lambda}$ where
$M_\lambda=\{w\in M|L(0)w=\lambda w\}$. Moreover, one requires that $M_\lambda$ is
finite dimensional and for fixed $\lambda\in\C$, $M_{\lambda+n}=0$
for sufficiently small integer $n$. It is known \cite{DLM3} that a $V$-module is always admissible.

\vskip.5cm
We now recall from \cite{DZ1} the definition of the Zhu's algebra $A(V)$ of a $\Z$-graded vertex operator superalgebra $V$. $A(V)$ is defined as $A(V):=V/O(V),$ where $O(V)={\rm span}\{u\circ v|u, v\in V\}$ and
$$u\circ v:=\Res_z \frac{(1+z)^{\wt u}}{z^2}Y(u, z)v$$ for homogeneous $u, v\in V$ and extended linearly. It is an associative algebra with multiplication defined as $$u*v:=\Res_z \frac{(1+z)^{\wt u}}{z}Y(u, z)v$$ for homogeneous $u, v\in V$ (see Theorem 3.1 of \cite{DZ1}).

For a weak $V$-module $M$, define the space $\Omega(M)$ as follows:
$$\Omega(M)=\{w\in M|u_{\wt u-1+n}w=0, u\in V, n>0\}.$$
Then the following results have been proved in Theorem 3.2 of \cite{DZ1}.
\begin{theorem}\label{zhu}
Let $V$ be a $\Z$-graded vertex operator superalgebra, and $M$ be a weak $V$-module. Then\\
(1) $\Omega(M)$ is an $A(V)$-module such that $v+O(V)$ acts as $v_{\wt v-1}$.\\
(2) If $M=\oplus_{n\in \Z_{\geq 0}}M(n)$ is an admissible $V$-module such that $M(0)\neq 0$, then $M(0)\subseteq \Omega(M)$ is an $A(V)$-submodule. Moreover, $M$ is irreducible if and only if $M(0)=\Omega(M)$ and $M(0)$ is a simple $A(V)$-module.\\
(3) The map $M\rightarrow M(0)$ gives a $1-1$ correspondence between the irreducible admissible $V$-modules and simple $A(V)$-modules.
\end{theorem}

We also need the following result, which can be proved by the similar argument as Lemma 4.3 of \cite{A}.
\begin{proposition}\label{finite1}
Let $V$ be a $C_2$-cofinite $\Z$-graded vertex operator superalgebra. Then $A(V)$ is finite dimensional.
\end{proposition}
\subsection{Affine vertex operator superalgebras}
In this subsection, we recall from \cite{FZ}, \cite{K2}, \cite{LL} some facts about  affine vertex operator superalgebras. Let $\g$ be a finite dimensional simple Lie superalgebra with a nondegenerate even supersymmetric invariant bilinear form $(\cdot|\cdot)$. The affine Lie superalgebra associated to $\g$ is defined on $\hat{\g}=\g\otimes \C[t^{-1}, t]\oplus \C K$ with Lie brackets
\begin{align*}
[x(m), y(n)]&=[x, y](m+n)+(x|y) m\delta_{m+n,0}K,\\
[K, \hat\g]&=0,
\end{align*}
for $x, y\in \g$ and $m,n \in \Z$, where $x(n)$ denotes $x\otimes t^n$.

For a complex number $k$, define the vacuum module of $\hat \g$ by
\begin{align*}
V^{k}(\g)=\Ind_{\g\otimes \C[t]\oplus \C K}^{\hat \g}\C_k,
\end{align*}
where $\C_k=\C\1$ is the $1$-dimensional $\g\otimes \C[t]\oplus \C K$-module such that $\g\otimes \C[t]$ acts as $0$ and $K$ acts as $k$.
\begin{theorem}[\cite{FZ,K2}]\label{affine}
Let $\g$ be a finite dimensional simple Lie superalgebra with a nondegenerate even supersymmetric invariant bilinear form $(\cdot|\cdot)$, $h^\vee$ be the dual Coxeter number of $(\g, (\cdot|\cdot))$ and $k$ be a complex number which is not equal to $-h^\vee$. Then $V^{k}(\g)$ is a vertex operator superalgebra such that the conformal vector $\omega$ is defined as follows:
$$\omega=\frac{1}{2(k+h^\vee)}\sum_ia^i(-1)b^i(-1)\1,$$where $\{a^i\}$ and $\{b^i\}$ are dual bases of $\g$ with respect to $(\cdot|\cdot)$. Moreover, the conformal weight of $x(-1)\1$ is $1$ for any $x\in \g$.
\end{theorem}

 Note that $V^{k}(\g)$ is a strong CFT type  $\Z$-graded vertex operator superalgebra. It is well-known that $V^{k}(\g)$ has a unique irreducible quotient module which is denoted by $L_{\g}(k, 0)$ (see \cite{K}).  Then $L_{\g}(k, 0)$ has a vertex operator superalgebra structure \cite{LL}. Moreover, $L_{\g}(k, 0)$ is also a $\Z$-graded vertex operator superalgebra.
\section{Integrability of $C_2$-cofinite vertex operator superalgebras }
\def\theequation{3.\arabic{equation}}
\setcounter{equation}{0}
Let $V$ be a vertex operator superalgebra of strong CFT type. It is known \cite{B} that the weight one subspace $V_1$ of $V$ has a Lie superalgebra structure defined by $[u, v]=u_0v$ for any $u, v\in V_1$. Moreover,  we have
\begin{lemma}
Let $V$ be a vertex operator superalgebra of strong CFT type. Then $V_1$ is a Lie superalgebra equipped with a supersymmetric invariant bilinear form $B(\cdot, \cdot)$ such that $B(u, v)\1=u_1v$ for any $u, v\in V_1$.
\end{lemma}
\pf We first prove that $B(\cdot, \cdot)$  is supersymmetric. By the  skew-symmetry property of vertex operator superalgebras (see (2.2.5) of \cite{Li}), we have $$Y(u, z)v=(-1)^{[ u][ v]}e^{zL(-1)}Y(v, -z)u,$$
for any $\Z_2$-homogeneous $u, v\in V$. It follows that $u_1v=(-1)^{[ u][ v]}v_1u$ for any $u, v\in V_1$. Hence, $B(u,v)=(-1)^{[ u][ v]}B(v, u)$ for any $u, v\in V_1$.
We next prove that $B(\cdot, \cdot)$  is invariant, i.e., $B([u, v], w)=B(u, [v, w])$. By (2.2.6) of \cite{Li}, we have for any  $\Z_2$-homogeneous $u, v, w\in V$,
\begin{align}\label{comm}
u_mv_nw-(-1)^{[ u][ v]}v_mu_nw=\sum_{i=0}^{\infty}\binom{m}{i}(u_iv)_{m+n-i}w.
\end{align}
Hence, we have for any  $\Z_2$-homogeneous $u, v, w\in V_1$, $$u_1v_0w-(-1)^{[ u][ v]}v_0u_1w=(u_0v)_{1}w+(u_1v)_{0}w.$$
This implies that $u_1v_0w=(u_0v)_{1}w$ for any  $\Z_2$-homogeneous $u, v, w\in V_1$. Then we have the invariant property.
\qed

\vskip.25cm
We now let $\g$ be a subalgebra of $V_1$. Suppose further that $\g$ is a simple Lie superalgebra equipped with a nondegenerate even supersymmetric invariant bilinear from $(\cdot|\cdot)$. Then there exists a complex number $k$ such that $B(u,v)=k(u|v)$ for any $u, v\in \g$. By the formula (\ref{comm}) we have for any $\Z_2$-homogeneous $u, v\in V_1$,
\begin{align*}
&u_mv_n-(-1)^{[ u][ v]}v_mu_n\\
&=(u_0v)_{m+n}+m(u_1v)_{m+n-1}\\
&=([u, v])_{m+n}+mB(u, v)\delta_{m+n, 0}\\
&=([u, v])_{m+n}+mk(u| v)\delta_{m+n, 0}.
\end{align*}
Thus, $V$ is a $\hat\g $-module of level $k$.

\vskip.5cm
We next assume further that $V$ is $C_2$-cofinite. Fix homogeneous $x^i\in V$ ($i$ ranging over an index set $I$) such that the coset $x^i+C_2(V), i\in I$, span a complement to $\1+C_2(V)$ in $V/C_2(V)$. Set $X=\{x^i|i\in I\}$. Then the following property of $V$ has been essentially proved in Proposition 8 of \cite{GN}.
\begin{theorem}\label{span}
Let $V$ be a vertex operator superalgebra which is of strong CFT type and $C_2$-cofinite. Then $V$ is spanned by $\1$ together with elements of the form $$x^1_{-n_1}\cdots x^k_{-n_k}\1$$ where $n_1>\cdots>n_k>0, x^i\in X$.
\end{theorem}
 As a consequence, by the same argument as in Theorem 3.1 of \cite{DM}, we have the following result, which was essentially proved in Theorem 6.1 of \cite{DH}.
  \begin{theorem}\label{keyDH}
  Let $V$ be a vertex operator superalgebra which is of strong CFT type and $C_2$-cofinite, $\g$ be a Lie subalgebra of $V_1$. Assume that $\g$ is a simple Lie algebra, then $V$ is an integrable $\hat \g$-module.
  \end{theorem}

  We are now ready to prove the main result in this section.
\begin{theorem}\label{inte}
Let $V$ be a vertex operator superalgebra which is of strong CFT type and $C_2$-cofinite, $\g$ be a  Lie subalgebra of $V_1$. Assume that $\g$ is a basic simple Lie superalgebra, then $\g$ is isomorphic to $osp(1|2n)$, $sl(1|n+1)$ or $osp(2|2n), n\geq 1$.
\end{theorem}
\pf Let $\g$ be a  basic simple Lie superalgebra. Assume that $\g$ is not isomorphic to $osp(1|2n)$, $sl(1|n+1)$ or $osp(2|2n), n\geq 1$. We choose the  nondegenerate even supersymmetric invariant bilinear from $(\cdot|\cdot)$ of $\g$ defined in Table 6.1 of \cite{KWaki}. Thus, there are simple Lie subalgebras $\f_1$ and $\f_2$ of $\g_{\bar 0}$ such that the restriction of $(\cdot|\cdot)$ on $\f_1$ is positive definite and the restriction of $(\cdot|\cdot)$ on $\f_2$ is negative definite.

 Consider the vertex subalgebra $\langle \g\rangle$ of $V$ generated by $\g$, then $\langle \g\rangle$ is isomorphic to a quotient module of  the vacuum module $V^k(\g)$ of $\hat \g$ for some complex number $k$.  By Theorem \ref{keyDH}, $V$ is an integrable module of $\hat \f_1$. Hence, $\langle \g\rangle$ is an integrable  module of $\hat \f_1$.  Note that the vacuum vector $\1$ is a highest weight vector for $\hat \f_1$. Let $W$ be the $\hat \f_1$-submodule of $\langle \g\rangle$  generated by $\1$. Then $W$ is an integrable highest weight module of  $\hat \f_1$. As a consequence, $W$ is an irreducible integrable highest weight module of  $\hat \f_1$. Therefore, by Lemma 10.1 of \cite{K}, $k$ must be a nonnegative number.

 On the other hand, $V$ is also an integrable module of $\hat \f_2$. Hence, $\langle \g\rangle$ is an integrable  module of $\hat \f_2$. Note that the restriction of $(\cdot|\cdot)$ on $\f_2$ is negative definite. By the similar argument as above, we can prove that $k$ must be a nonpositive number. As a consequence, $k=0$. This forces that $W$ viewed as an $\hat \f_1$-module is isomorphic to  $\C$. However, $W$ contains a subspace isomorphic to $\f_1$,  this is a contradiction. Thus, $\g$ must be isomorphic to $osp(1|2n)$, $sl(1|n+1)$ or $osp(2|2n), n\geq 1$.
\qed
\section{Irreducible modules of vertex operator superalgebras $L_{sl(1| n)}(k, 0)$ and $L_{osp(2| 2n)}(k, 0)$}
\def\theequation{4.\arabic{equation}}
\setcounter{equation}{0}
\subsection{Irreducible modules of vertex operator superalgebra $L_{sl(1| n)}(k, 0)$} In this subsection, for any positive integer $k$, we will show that $L_{sl(1| n)}(k, 0)$ has infinitely many nonisomorphic irreducible admissible modules. As a result, we will show that $L_{sl(1| n)}(k, 0)$ is not $C_2$-cofinite.

Consider the simple Lie superalgebra $sl(1|n)$, $n\geq 2$, which consists of block matrices of the form $X=\left(\begin{array}{cc}A & B\\ C & D\end{array}\right)$ such that $A$ is a $(1\times 1)$-matrix, $B$ is a $(1\times n)$-matrix, $C$ is a $(n\times 1)$-matrix, $D$ is a $(n\times n)$-matrix and that ${\rm Str}(X):=\tr(A)-\tr(D)=0$ (cf. \cite{M}). For any $X, Y\in sl(1|n)$, define $(X|Y)=-{\rm Str}(XY)$, then $(\cdot| \cdot)$ is a nondegenerate even supersymmetric invariant bilinear form of $sl(1|n)$. We use $I$ to denote the identity matrix and let $h$ be the matrix  $\left(\begin{array}{cc}1 & 0\\ 0 & \frac{1}{n}I\end{array}\right)$. Then $h$ is an element of $sl(1|n)_{\bar 0}$. Let $\f$ be the set of matrices of the form $\left(\begin{array}{cc}0 & 0\\ 0 & D\end{array}\right)$ such that $\tr D=0$, then $\f$ is a subalgebra of $sl(1|n)_{\bar 0}$. Moreover, $\f$ is isomorphic to the simple Lie algebra $sl_n$. Thus the even part of $sl(1|n)$ is $sl_n\oplus \C h$. Let $\h$ be the subset of $sl(1|n)$ consisting of diagonal matrices. Then $\h$ is a Cartan subalgebra of $sl(1|n)$. We next describe the root system of $sl(1|n)$. Let $\epsilon_i,0\leq i \leq n,$ be the linear functionals on $\h$ whose values on the diagonal matrix $$a={\rm diag}(a_0,\cdots, a_{n})$$ are given by $\epsilon_i(a)=a_i$. Then the root system of $sl_n$ is $\{\pm(\epsilon_i-\epsilon_j)|1\leq i<j\leq n\}$ and the set of odd roots of $sl(1|n)$ is $\{\pm(\epsilon_0-\epsilon_i)|1\leq i\leq n\}$ (see subsection 2.2 of \cite{M}). Moreover, the bilinear form $(\cdot| \cdot)$ is determined by $(\epsilon_i| \epsilon_j)=0$ if $i\neq j$, $(\epsilon_0| \epsilon_0)=-1$ and  $(\epsilon_i| \epsilon_i)=1$ for $1\leq i\leq n$.  In the following, we choose the following subset of simple roots of $sl(1|n)$:
$$\Pi=\{\epsilon_0-\epsilon_1,\epsilon_1-\epsilon_2, \cdots, \epsilon_{n-1}-\epsilon_n\}.$$ Note that the highest root of $sl_n$ is  $\theta=\epsilon_1-\epsilon_n$. Let $e_\theta$ be a highest root vector of $sl_n$. Then we have the following result, which was proved in Theorem 5.4.1 and Corollary 5.4.3 of \cite{GS}.
\begin{theorem}\label{ideal}
Let $k$ be a positive integer. Then $L_{sl(1| n)}(k, 0)=V^k(sl(1| n))/I$, where $I$ is the submodule of $V^k(sl(1| n))$ generated by $e_{\theta}(-1)^{k+1}\1$.
\end{theorem}

Thus, to find irreducible admissible modules of $L_{sl(1| n)}(k, 0)$, it is equivalent to find irreducible admissible modules of $V^k(sl(1| n))$ such that $Y(e_{\theta}(-1)^{k+1}\1, z)=0$. We next show that $L_{sl(1| n)}(k, 0)$ has infinitely many irreducible admissible modules. Consider the affine Kac-Moody superalgebra $sl(1| n)^{(1)}:=\hat{sl(1| n)}\oplus \C d$. Set $\tilde{\h}=\h\oplus \C K\oplus \C d$, then $\tilde{\h}$ is a Cartan subalgebra of $sl(1| n)^{(1)}$. We use $\delta$ to denote the smallest positive imaginary root of $sl(1| n)^{(1)}$, and choose the following subset $\Sigma$ of simple roots of $sl(1| n)^{(1)}$:
 $$\Sigma=\{\delta-\epsilon_0+\epsilon_n, \epsilon_0-\epsilon_1,\epsilon_1-\epsilon_2, \cdots, \epsilon_{n-1}-\epsilon_n\}.$$
For $\lambda\in (\tilde{\h})^*$, let $L_{\Sigma}(\lambda)$ be the irreducible highest weight module of $sl(1| n)^{(1)}$ of highest weight $\lambda$.  Recall that $L_{\Sigma}(\lambda)$ is an integrable $sl(1| n)^{(1)}$-module if $L_{\Sigma}(\lambda)$ is integrable over $sl_n^{(1)}$ and locally finite over the Cartan subalgebra $\h$ (see \cite{GS}, \cite{KWaki1}). Set
 $$b_0=(\lambda| \delta-\epsilon_0+\epsilon_n),~b_i=(\lambda| \epsilon_{i-1}-\epsilon_i),~ 1\leq i\leq n.$$ Then the following result was proved in \cite{KWaki1} (see also subsection 2.3.1  of \cite{GS}).
\begin{theorem}\label{inte1}
Let $k$ be a positive integer. Then  $L_{\Sigma}(\lambda)$ is an integrable $sl(1| n)^{(1)}$-module of level $k$ if and only if \\
(i) $b_i\in \Z_{\geq 0}$ for $i\geq 2$.\\
(ii) $b_0+b_1\in \Z_{>0}$ or $b_0=b_1=0$.\\
(iii)  $b_0+b_1+\cdots+b_n=k$.
\end{theorem}

Furthermore, we have
\begin{proposition}\label{irr}
Let $k$ be a positive integer, $L_{\Sigma}(\lambda)$ be an integrable $sl(1| n)^{(1)}$-module of level $k$. Then $L_{\Sigma}(\lambda)$ is an irreducible admissible $L_{sl(1| n)}(k, 0)$-module.
\end{proposition}
\pf First, $L_{\Sigma}(\lambda)$ is a weak $V^k(sl(1| n))$-module (see subsection 5.3 of \cite{GS}). We next show that $Y(e_{\theta}(-1)^{k+1}\1, z)=0$ on  $L_{\Sigma}(\lambda)$. Since $L_{\Sigma}(\lambda)$ is an integrable $sl(1| n)^{(1)}$-module of level $k$, it is an integrable $\hat{sl}_ n$-module of level $k$. By Remark 3.9 of \cite{DLM}, $L_{\Sigma}(\lambda)$ viewed as an $\hat{sl}_ n$-module is a direct sum of irreducible highest weight integrable $\hat{sl}_ n$-modules. By Proposition 6.6.21 of \cite{LL}, we have $Y(e_{\theta}(-1)^{k+1}\1, z)=0$ on  $L_{\Sigma}(\lambda)$. As a consequence, $L_{\Sigma}(\lambda)$ is a weak $L_{sl(1| n)}(k, 0)$-module by Theorem \ref{ideal}.

We next show that $L_{\Sigma}(\lambda)$ is an admissible $L_{sl(1| n)}(k, 0)$-module. Note that $d$ acts diagonally on $L_{\Sigma}(\lambda)$ with eigenvalues bounded from above. By the discussion in Subsection 5.6 of \cite{GS}, $L_{\Sigma}(\lambda)$ is an admissible $V^k(sl(1|n))$-module (see also Remark 6.6.1 of \cite{LL}). This implies that $L_{\Sigma}(\lambda)$ is an admissible $L_{sl(1| n)}(k, 0)$-module.
\qed

\vskip.25cm
Note that by Theorem \ref{inte1} there are infinitely many integrable $sl(1| n)^{(1)}$-modules of level $k$. Thus, by Proposition  \ref{irr}, there are infinitely many irreducible admissible $L_{sl(1| n)}(k, 0)$-modules. By Proposition \ref{finite1}, we have
\begin{theorem}\label{non1}
Let $k$ be a positive integer. Then the vertex operator superalgebra $L_{sl(1| n)}(k, 0)$ is not $C_2$-cofinite.
\end{theorem}
\subsection{Irreducible modules of vertex operator superalgebra  $L_{osp(2| 2n)}(k, 0)$} In this subsection, for a positive integer $k$, we will show that $L_{osp(2| 2n)}(k, 0)$ has infinitely many nonisomorphic irreducible admissible modules. As a result, we will show that $L_{osp(2| 2n)}(k, 0)$ is not $C_2$-cofinite.

Consider the simple Lie superalgebra $osp(2| 2n)$, $n\geq 1$, which consists of block matrices of the form $$X=\left(\begin{array}{cccc}a &0&y& y_1\\ 0&-a&z &z_1 \\-z_1^t&-y_1^t&d&e\\z^t&y^t&f&-d^t\end{array}\right)$$ such that $a$ is a complex number, $y,y_1,z,z_1$ are $(1\times n)$-matrices, $e,f$ are symmetric $(n\times n)$-matrices and $d$ is a $(n\times n)$-matrix (cf. \cite{M}). For any $X, Y\in osp(2| 2n)$, define $(X|Y)=-{\rm Str}(XY)$, then $(\cdot| \cdot)$ is a nondegenerate even supersymmetric invariant bilinear form of $osp(2| 2n)$. Set   $$h=\left(\begin{array}{cccc}1 &0&0& 0\\ 0&-1&0 &0 \\0&0&0&0\\0&0&0&0\end{array}\right).$$ Then $h$ is an element of $osp(2| 2n)_{\bar 0}$. Let $\f$ be the set of matrices of the form $$\left(\begin{array}{cccc}0&0&0& 0\\ 0&0&0 &0 \\0&0&d&e\\0&0&f&-d^t\end{array}\right)$$ such that $e,f$ are symmetric $(n\times n)$-matrices and $d$ is a $(n\times n)$-matrix.  Then $\f$ is a subalgebra of $osp(2| 2n)_{\bar 0}$. Moreover, $\f$ is isomorphic to the simple Lie algebra of type $C_n$. Thus the even part of $osp(2| 2n)$ is $\f\oplus \C h$. Let $\h$ be the subset of $osp(2| 2n)$ consisting of diagonal matrices. Then $\h$ is a Cartan subalgebra of $osp(2| 2n)$. We next describe the root system of $osp(2| 2n)$. Let $\epsilon_i,0\leq i \leq n,$ be the linear functionals on $\h$ whose values on the diagonal matrix $$a={\rm diag}(a_1,\cdots, a_{2n+2})$$ are given by $\epsilon_i(a)=a_i$. Then the root system of $\f$ is $\{\pm(\epsilon_i\pm\epsilon_j), \pm2\epsilon_i|3\leq i<j\leq n+2\}$ and the set of odd roots of $osp(2| 2n)$ is $\{\pm(\epsilon_1-\epsilon_i)|3\leq i\leq n+2\}$ (see Subsection 2.3.3 of \cite{M}). Moreover, the bilinear form $(\cdot| \cdot)$ is determined by $(\epsilon_i| \epsilon_j)=0$ if $i\neq j$, $(\epsilon_1| \epsilon_1)=-\frac{1}{2}$ and  $(\epsilon_i| \epsilon_i)=\frac{1}{2}$ for $3\leq i\leq n+2$.  In the following, we choose the following subset of simple roots of $osp(2| 2n)$:
$$\Pi=\{\epsilon_1-\epsilon_3,\epsilon_3-\epsilon_4, \cdots, \epsilon_{n+1}-\epsilon_{n+2}, 2\epsilon_{n+2}\}.$$ Note that the highest root of $\f$ is  $\theta=2\epsilon_3$. Let $e_\theta$ be a highest root vector of $\f$. Then we have the following result, which was proved in Theorem 5.4.1 and Corollary 5.4.3 of \cite{GS}.
\begin{theorem}\label{ideal1}
Let $k$ be a positive integer. Then $L_{osp(2| 2n)}(k, 0)=V^k(osp(2| 2n))/I$, where $I$ is the submodule of $V^k(osp(2| 2n))$ generated by $e_{\theta}(-1)^{k+1}\1$.
\end{theorem}

Thus, to find irreducible admissible modules of $L_{osp(2| 2n)}(k, 0)$, it is equivalent to find irreducible admissible modules of $V^k(osp(2| 2n))$ such that $Y(e_{\theta}(-1)^{k+1}\1, z)=0$. We next show that $L_{osp(2| 2n)}(k, 0)$ has infinitely many irreducible admissible modules. Consider the affine Kac-Moody superalgebra $osp(2| 2n)^{(1)}:=\hat{osp(2| 2n)}\oplus \C d$. Set $\tilde{\h}=\h\oplus \C K\oplus \C d$, then $\tilde{\h}$ is a Cartan subalgebra of $osp(2| 2n)^{(1)}$. We use $\delta$ to denote the smallest positive imaginary root of $osp(2| 2n)^{(1)}$, and choose the following subset $\Sigma$ of simple roots of $osp(2| 2n)^{(1)}$:
 $$\Sigma=\{\delta-\epsilon_1-\epsilon_3, \epsilon_1-\epsilon_3,\epsilon_3-\epsilon_4, \cdots, \epsilon_{n+1}-\epsilon_{n+2}, 2\epsilon_{n+2}\}.$$
For $\lambda\in (\tilde{\h})^*$, let $L_{\Sigma}(\lambda)$ be the irreducible highest weight module of $osp(2| 2n)^{(1)}$ of highest weight $\lambda$.  Recall that $L_{\Sigma}(\lambda)$ is an integrable $osp(2| 2n)^{(1)}$-module if $L_{\Sigma}(\lambda)$ is integrable over $\f^{(1)}$ and locally finite over the Cartan subalgebra $\h$ (see \cite{GS}, \cite{KWaki1}). Set
 $$b_0=(\lambda| \delta-\epsilon_1-\epsilon_3),~b_1=(\lambda| \epsilon_1-\epsilon_3),~b_i=2(\lambda| \epsilon_{i-1}-\epsilon_i),~ 4\leq i\leq n+2, ~b_{n+3}=(\lambda| 2\epsilon_{n+2}).$$ Then the following result was proved in Theorem 8.1 of \cite{KWaki1}.
\begin{theorem}\label{inte2}
Let $k$ be a positive integer. Then  $L_{\Sigma}(\lambda)$ is an integrable $osp(2| 2n)^{(1)}$-module of level $k$ if and only if \\
(i) $b_i\in \Z_{\geq 0}$ for $i\geq 4$.\\
(ii) $b_0+b_1\in \Z_{>0}$ or $b_0=b_1=0$.\\
(iii)  $b_0+b_1+b_4+\cdots+b_{n+2}+b_{n+3}=k$.
\end{theorem}

Furthermore, by the similar argument as that in Proposition  \ref{irr}, we have the following result.
\begin{proposition}\label{irr1}
Let $k$ be a positive integer, $L_{\Sigma}(\lambda)$ be an integrable $osp(2| 2n)^{(1)}$-module of level $k$. Then $L_{\Sigma}(\lambda)$ is an irreducible admissible $L_{osp(2| 2n)}(k, 0)$-module.
\end{proposition}

Note that by Theorem \ref{inte2} there are infinitely many integrable $osp(2| 2n)^{(1)}$-modules of level $k$. Thus, by Proposition \ref{irr1}, there are infinitely many irreducible admissible $L_{osp(2| 2n)}(k, 0)$-modules. By Proposition \ref{finite1}, we have
\begin{theorem}\label{non2}
Let $k$ be a positive integer. Then the vertex operator superalgebra $L_{osp(2| 2n)}(k, 0)$ is not $C_2$-cofinite.
\end{theorem}

\section{$C_2$-cofiniteness of affine vertex operator superalgebras}
\def\theequation{5.\arabic{equation}}
\setcounter{equation}{0}
In this section, for any positive integer $n$, we will show that the vertex operator superalgebra $L_{osp(1| 2n)}(k, 0)$ is  $C_2$-cofinite if $k$ is a positive integer. This result has been known for many years (see Section 0.4 of  \cite{GK}).  For completeness, we give a proof here.

 We need the following result, which can be proved by the similar argument as  in Lemma 3.8 of \cite{DLM2}.
\begin{lemma}\label{basic1}
Let $V$ be a vertex operator superalgebra. Then $C_2(V)$ is closed under the operators $v_0$ and $v_{-1}$ for any $v\in V$.
\end{lemma}

We also need the following fact which was proved in Lemma 12.1 of \cite{DLM2}. \begin{lemma}\label{basic2}
Let $V$ be a vertex operator superalgebra. Then $C_2(V)$ contains $u_{-m}V$ for any $ u\in V$ and $m\geq 2 $.
\end{lemma}

We now let $\g$ be a simple Lie superalgebra which is isomorphic to $osp(1| 2n)$, $n\geq 1$. Let $\g_{\bar 0}$ be  the even part of $\g$, then $\g_{\bar 0}$ is isomorphic to the simple Lie algebra $sp(2n)$. In the following, we fix a Cartan subalgebra $\h$ of $\g_{\bar 0}$. Note that there is a nondegenerate supersymmetric invariant bilinear form $(\cdot|\cdot)$ of $\g$ such that the restriction of $(\cdot|\cdot)$ on $\g_{\bar 0}$ is the normalized invariant bilinear form of $\g_{\bar 0}$. Then we have
\begin{proposition}
Let $k$ be a positive integer. Then for any root vector  $e_\alpha$ of $\g_{\bar 0}$,  $Y(e_\alpha, z)^{tk+1}=0$ acting on  $L_{\g}(k, 0)$, where $t=1$ if $\alpha$ is a long root of $\g_{\bar 0}$ and $t=2$ if $\alpha$ is a short root of $\g_{\bar 0}$.
\end{proposition}
\pf The argument of the proof is similar to that in Proposition 5.2.1 of \cite{Li}. To prove the statement, it is enough to prove $e_\alpha^{tk+1}\1=0$. From the standard semisimple Lie algebra theory (cf. \cite{H}) we can embed $sl_2$ into $\g_{\bar 0}$ as $\g^\alpha$ linearly spanned by $e_\alpha, f_\alpha, h_\alpha$.  It is known  \cite{GS}, \cite{K1} that  $L_{\g}(k, 0)$ is an integrable $\hat{\g}$-module, then  $L_{\g}(k, 0)$ is an integrable $\hat{\g^\alpha}$-module. As a result, $\1$ generates an integrable $\hat{\g^\alpha}$-module $W$. In particular, $W$ is an irreducible $\hat{\g^\alpha}$-module. Set $t=\frac{2}{(\alpha, \alpha)}$, then $(h_\alpha, h_\alpha)=\frac{4}{(\alpha, \alpha)}=2t$. As a result, we have $[(f_\alpha)_1, (e_\alpha)_{-1}]=h_\alpha+tk$. Hence, $(f_\alpha)_1(e_\alpha)_{-1}^{tk+1}\1=0$. Note also that $(e_\alpha)_0(e_\alpha)_{-1}^{tk+1}\1=0$. Hence, if $(e_\alpha)_{-1}^{tk+1}\1\neq 0$, it generates a $\hat{\g^\alpha}$-submodule of $W$, this is a contradiction. Then $(e_\alpha)_{-1}^{tk+1}\1= 0$.
\qed

Furthermore, by the similar argument as in Lemma 3.6 of \cite{DLM}, we have
\begin{corollary}\label{nilp}
Let  $k$ be a positive integer. Then there is a basis of $\{a^1, \cdots, a^m\}$ of $\g_{\bar 0}$ such that for $1\leq i\leq m$ we have $$[Y(a^i, z), Y(a^i, z)]=0~{\rm and}~Y(a^i, z)^{3k+1}=0$$ as operators on $L_{\g}(k, 0)$.
\end{corollary}

We are now ready to prove the following result.
\begin{theorem}
Let $k$ be a positive integer. Then  the vertex operator superalgebra $L_{osp(1| 2n)}(k, 0)$ is $C_2$-cofinite.
\end{theorem}
\pf The argument of the proof is similar to that in Proposition 12.6 of \cite{DLM2}. By definition,  $$V^k(\g)=U(\hat \g)\otimes_{U(\sum_{m=0}^{\infty} t^m\otimes \g\oplus \C K)}\C\cong U(\sum_{m=1}^{\infty }t^{-m}\otimes \g)~{\rm (linearly)}.$$ By Lemma \ref{basic2}, $C_2(L_\g(k, 0))$ contains $a_{-m}L_\g(k, 0)$ for any $a\in \g$ and $m\geq 2$. Thus, $L_\g(k, 0)=C_2(L_\g(k, 0))+U(t^{-1}\otimes \g)\1$. It is enough to show that $C_2(L_\g(k, 0))$ contains $$x_1^{n_1}(-1)\cdots x_s^{n_s}(-1)\1,$$ whenever $n_i\geq 0$ and $n_1+\cdots +n_s$ is large enough; here $x_1, \cdots, x_s$ is a basis of $\g$.

First, let $a$ be a root vector in the odd part of $\g$. Then we have $a(-1)^2=\frac{1}{2}[a, a](-2)$. As a consequence,  $a(-1)^2b$ is contained in $C_2(L_\g(k, 0))$ for any $b\in L_\g(k, 0)$.  On the other hand, by Corollary \ref{nilp}, there is a basis $\{a^1, \cdots, a^m\}$ of $\g_{\bar 0}$ such that for $1\leq i\leq m$ we have $$[Y(a^i, z), Y(a^i, z)]=0~{\rm and}~Y(a^i, z)^{3k+1}=0$$ as operators on $L_{\g}(k, 0)$. The constant term of $Y(a^i,z)^{3k+1}\1$ is equal to $a^i(-1)^{3k+1}\1+r$, where $r$ is a sum of products of the form $a^i(n_1)^{e_1}\cdots a^i(n_{3k+1})^{e_{3k+1}}1$ with some $n_j\leq -2$. Since the operators $a^i(n_j)$ commute, we have $r\in C_2(L_\g(k, 0))$. Hence, $a^i(-1)^{3k+1}\1\in C_2(L_\g(k, 0))$. By Lemmas \ref{basic1}, \ref{basic2}, we can conclude that $C_2(L_\g(k, 0))$ contains $x_1^{n_1}(-1)\cdots x_s^{n_s}(-1)\1$ whenever $n_i\geq 3k+1$ for some $i$. This completes the proof.
\qed

Combining with Theorems \ref{inte}, \ref{non1}, \ref{non2}, we immediately obtain the following result, which has also been known for many years (see Section 0.4 of  \cite{GK}).
\begin{theorem}\label{main}
Let $\g$ be a basic simple Lie superalgebra. Then $L_{\g}(k, 0)$ is $C_2$-cofinite if and only if $\g$ is isomorphic to $osp(1|2n)$ and $k$ is a nonnegative integer.
\end{theorem}

\section{Category of $L_{G(3)}(1, 0)$-modules}
\def\theequation{6.\arabic{equation}}
\setcounter{equation}{0} 

In this section, we give a proof of the fact that the category of $L_{G(3)}(1, 0)$-modules is semisimple. As a consequence, we show that $L_{G(3)}(1,0)$ is a vertex operator superalgebra such that the category of $L_{G(3)}(1,0)$-modules is semisimple but $L_{G(3)}(1,0)$ is not $C_2$-cofinite.

Let $\g_{\bar 0}=sl_2\oplus G_2$ and $\g_{\bar 1}=V_1\otimes V$, where $V_1$ denotes the standard two-dimensional module of $sl_2$ and $V$ denotes the unique $7$-dimensional irreducible $G_2$-module. Then there is a simple Lie superalgebra structure on $G(3):=\g_{\bar 0}\oplus \g_{\bar 1}$ (see Subsection 1.4.2 of  \cite{S}). We next describe the root system of $G(3)$. First, recall \cite{C} that $G_2$ has a Cartan subalgebra $\h_2$ such that we may choose the following positive root system of $G_2$:
$$\Delta_2^+=\{\alpha_1, \alpha_2, \alpha_1+\alpha_2, \alpha_1+2\alpha_2, \alpha_1+3\alpha_2, 2\alpha_1+3\alpha_2\}.$$ The bilinear form $(\cdot|\cdot)$ on $\h_2$ is determined by  $(\alpha_1|\alpha_2)=-1,$ $(\alpha_2|\alpha_1)=-1$, $(\alpha_1|\alpha_1)=2$ and $(\alpha_2| \alpha_2)=\frac{2}{3}$. The set of weights of $G_2$-module $V$  is as follows:
$$\{\alpha_1+2\alpha_2, \alpha_1+\alpha_2, \alpha_2, 0, -\alpha_2, -\alpha_1-\alpha_2, -\alpha_1-2\alpha_2\}.$$
We next choose the Cartan subalgebra $\h_1$ of $sl_2$ such that the root system of $sl_2$ is $\{\alpha_0, -\alpha_0\}$. We now choose the following  positive root system of $sl_2$: $$\Delta_1^+=\{\alpha_0\}.$$ The bilinear form $(\cdot|\cdot)$  on $\h_1$ is determined by $(\alpha_0|\alpha_0)=-\frac{8}{3}$. The set of weights of $sl_2$-module $V_1$ is as follows:
$$\{\frac{1}{2}\alpha_0, -\frac{1}{2}\alpha_0\}.$$
Note that  $\h_1\oplus \h_2$ is a Cartan subalgebra of $G(3)$. Then we may choose the following positive root system of $G(3)$:
$$\Delta^+=\Delta^+_1\cup \Delta^+_0,$$
where $\Delta^+_1=\{\frac{1}{2}\alpha_0, \frac{1}{2}\alpha_0\pm\alpha_2, \frac{1}{2}\alpha_0\pm(\alpha_1+2\alpha_2), \frac{1}{2}\alpha_0\pm(\alpha_1+\alpha_2)\}$ and
$\Delta^+_0=\{\alpha_0, \alpha_1, \alpha_2, \alpha_1+\alpha_2, \alpha_1+2\alpha_2, \alpha_1+3\alpha_2, 2\alpha_1+3\alpha_2\}$ (see Subsection 10.9 of \cite{GK}).

Note that $\theta_2=2\alpha_1+3\alpha_2$ is the highest root of $G_2$ and that the bilinear form $(\cdot|\cdot)$ on $G(3)$ is normalized such that $(\theta_2|\theta_2)=2$. We now let $e_{\theta_2}$ be a highest root vector $G_2$. Then the following result has been proved in Theorem 5.4.1 and Corollary 5.4.3 of \cite{GS}.
\begin{theorem}
The vertex operator superalgebra $L_{G(3)}(1, 0)$ is isomorphic to $V^1(G(3))/I$, where $I$ is the submodule generated by $e_{\theta_2}(-1)^2\1$.
\end{theorem}

The Zhu's algebra $A(L_{G(3)}(1, 0))$  of the vertex operator superalgebra $L_{G(3)}(1, 0)$ has been determined in Subsection 5.6.2 of \cite{GS}.
\begin{proposition}\label{zhu2}
The algebra $A(L_{G(3)}(1, 0))$  is isomorphic to $U(G(3))/(e_{\theta_2}^2)$.
\end{proposition}

We now classify irreducible $L_{G(3)}(1, 0)$-modules.
\begin{proposition}\label{classify}
$L_{G(3)}(1, 0)$ is the only irreducible $L_{G(3)}(1, 0)$-module.
\end{proposition}
\pf Let $M$ be an irreducible $L_{G(3)}(1, 0)$-module. Then $M$ may be viewed as an irreducible $\hat{G(3)}$-module of level $1$. Since $M$ is an irreducible $L_{G(3)}(1, 0)$-module, $M$  has the decomposition $$M=\bigoplus_{\lambda\in\C}
M_{\lambda},$$ where
$M_\lambda=\{w\in M|L(0)w=\lambda w\}$ and $\dim M_{\lambda}<\infty$ for any $\lambda\in \C$. Let $\lambda_0$ be a number such that $M_{\lambda_0}\neq 0$ but $M_{\lambda_0-n}=0$ for $n\in \Z_{> 0}$. By Theorem \ref{zhu}, $M_{\lambda_0}$ is an $A(L_{G(3)}(1, 0))$-module. Thus, by Proposition \ref{zhu2}, $M_{\lambda_0}$ may be viewed as a finite dimensional $G(3)$-module.  This implies that $M$ is an irreducible highest weight module of $\hat{G(3)}$.

We next show that $M$ is an integrable $\hat{G(3)}$-module, that is, $M$ is an integrable $\hat{G_2}$-module and locally finite with respect to $sl_2$ (cf. \cite{KWaki}). By Theorem 5.3.1 of \cite{GS}, $M$ is an integrable $\hat{G_2}$-module. Note that $M_\lambda$ is closed under the action of $G(3)$. Hence, the action of $G(3)$ on $M$ is locally finite. Thus, $M$ is an integrable irreducible highest weight  $\hat{G(3)}$-module. On the other hand, by Remark 6.3 of \cite{KWaki}, $L_{G(3)}(1, 0)$ is the only integrable  irreducible highest weight $\hat{G(3)}$-module of level $1$. Thus, $M$ is isomorphic to $L_{G(3)}(1, 0)$.
\qed

\vskip.25cm
 As a corollary, we have the following result.
\begin{corollary}\label{dimen}
 Let $M=\oplus_{n\in \Z_{\geq 0}}M_{\lambda+n}$ be an $L_{G(3)}(1, 0)$-module, where $M_{\lambda+n}=\{w\in M|L(0)w=(\lambda+n) w\}$ and $\lambda\in \C$. Then $\lambda=0$ and $M_\lambda$ contains an $A(L_{G(3)}(1, 0))$-submodule isomorphic to $\C$.
\end{corollary}
\pf Since $M$ is an $L_{G(3)}(1, 0)$-module, $M_\lambda$ is a finite dimensional $A(L_{G(3)}(1, 0))$-module by Theorem \ref{zhu}. Then there exists an irreducible $A(L_{G(3)}(1, 0))$-submodule $N$ of $M_\lambda$. Let $W$ be the $L_{G(3)}(1, 0)$-submodule of $M$ generated by $N$. Then $W$ is also an $L_{G(3)}(1, 0)$-module. Let $U$ be the sum of all the $L_{G(3)}(1, 0)$-submodules of $W$ which have zero intersection with $N$. It follows that $W/U$ is an irreducible $L_{G(3)}(1, 0)$-module. By Proposition \ref{classify}, we have $W/U$ is isomorphic to  $L_{G(3)}(1, 0)$.  For $\mu\in \C$, set  $(W/U)_{\mu}=\{w\in W/U|L(0)w=\mu w\}$. Note that $(W/U)_{\lambda}\cong N$. Therefore, $N$ viewed as an $A(L_{G(3)}(1, 0))$-module isomorphic to $\C$. It follows from Theorem \ref{zhu} that $\lambda=0$.
\qed

\vskip.25cm
To prove that the category of $L_{G(3)}(1, 0)$-modules is semisimple, we need the following result (see Lemma 2.3 of \cite{Ab}).
\begin{lemma}\label{weight}
Let $V$ be a $\Z$-graded vertex operator superalgebra, $v$ be a homogeneous vector in $V$, and $0\to M^1\to M\to M^2\to 0$ be a short exact sequence of weak $V$-modules. Suppose that $v_{\wt v-1}$ acts semisimply on both $M^1$ and $M^2$. Then $M$ is a direct sum of generalized eigenspaces for $v_{\wt v-1}$. Furthermore, the eigenspace of $M^2$ of eigenvalue $h$ is the image of the generalized eigenspace of $M$ of eigenvalue $h$.
\end{lemma}

We also need the following result.
\begin{lemma}
Let $0\to L_{G(3)}(1, 0)\to M\to L_{G(3)}(1, 0)\to 0$ be a short exact sequence of weak $L_{G(3)}(1, 0)$-modules. Then $M$ viewed as an $L_{G(3)}(1, 0)$-module is isomorphic to $L_{G(3)}(1, 0)\oplus L_{G(3)}(1, 0)$.
\end{lemma}
\pf By Lemma \ref{weight}, $M$ has the following decomposition $$M=\bigoplus_{\lambda\in\Z_{\geq 0}}
M_{\lambda},$$ where
$M_\lambda=\{w\in M|(L(0)-\lambda)^kw=0~{\rm for~ some }~k\}$ and $\dim M_{\lambda}<\infty$ for any $\lambda\in \Z_{\geq 0}$. By Theorem \ref{zhu}, we have $M_0$ is a two-dimensional $A(L_{G(3)}(1, 0))$-module.

We next show that $M_0$ viewed as an $A(L_{G(3)}(1, 0))$-module has the decomposition $M_0=\C\oplus \C$. Let $v$ be the image of the vacuum vector $\1$ of $L_{G(3)}(1, 0)$ in $M$, and $w$ be a vector of $M_0$ such that $M_0=\C v\oplus \C w$. By Proposition \ref{zhu2}, $M_0$ may be viewed as a finite dimensional $G(3)$-module. Therefore, $M_0$ viewed as a $G(3)_{\bar 0}=sl_2\oplus G_2$-module is completely reducible. Note that $\C v$ is an $A(L_{G(3)}(1, 0))$-submodule of $M_0$. This forces that  the action of $G(3)_{\bar 0}=sl_2\oplus G_2$ on $M_0$ must be trivial. In particular, $\h_1\oplus \h_2$ acts trivially on $M_0$.  As a consequence, $G(3)_{\bar 1}=V_1\otimes V$ must act trivially on $M_0$. Thus, $M_0$ viewed as an $A(L_{G(3)}(1, 0))$-module has the decomposition $M_0=\C\oplus \C$.

We now let $W$ be the $L_{G(3)}(1, 0)$-submodule of $M$ generated by $w$. Then $W$ also has the following decomposition $$W=\bigoplus_{\lambda\in\Z_{\geq 0}}
W_{\lambda},$$ where
$W_\lambda=\{w\in W|(L(0)-\lambda)^kw=0~{\rm for~ some }~k\}$. Since $\C w$ is an $A(L_{G(3)}(1, 0))$-submodule of $M_0$, we have $W_0=\C w$ by Proposition 4.1 of \cite{DM1}. Since $N$ is irreducible, we have $N\cap W=0$. This implies that $W$ is isomorphic to $L_{G(3)}(1, 0)$. Thus, $M$ viewed as an $L_{G(3)}(1, 0)$-module is isomorphic to $L_{G(3)}(1, 0)\oplus L_{G(3)}(1, 0)$.
 \qed

\vskip.25cm
As a consequence, by the similar argument as that in Lemma 2.6 of \cite{Ab}, we have the following
\begin{lemma}\label{split2}
Let $0\to \oplus_{i\in I} N^i \to M\to L_{G(3)}(1, 0)\to 0$ be a short exact sequence of weak $L_{G(3)}(1, 0)$-modules. Suppose that each $N^i$ is isomorphic to  $L_{G(3)}(1, 0)$. Then $M$ viewed as an $L_{G(3)}(1, 0)$-module is isomorphic to $(\oplus_{i\in I} N^i)\oplus L_{G(3)}(1, 0)$.
\end{lemma}

We are now ready to prove the main result in this section.
\begin{theorem}\label{semisimple}
The category of $L_{G(3)}(1, 0)$-modules is semisimple.
\end{theorem}
\pf Let $M$ be an $L_{G(3)}(1, 0)$-module. Then $M$ has the following decomposition $$M=\bigoplus_{\lambda\in\C}
M_{\lambda},$$ where
$M_\lambda=\{w\in M|(L(0)-\lambda)w=0\}$ and $\dim M_{\lambda}<\infty$ for any $\lambda\in \C$. Let $N$ be the sum of all irreducible $L_{G(3)}(1, 0)$-submodules of  $M$. If $M\neq N$, let $\lambda_0$ be a number such that $M_{\lambda_0}\neq N_{\lambda_0}$ but $M_{\lambda_0-n}= N_{\lambda_0-n}$ for $n\in \Z_{> 0}$. Then $M/N$ is an $L_{G(3)}(1, 0)$-module. For $\lambda\in \C$, set $(M/N)_\lambda=\{w\in M/N|(L(0)-\lambda)w=0\}$. Then we have $(M/N)_{\lambda_0}\neq 0$ and $(M/N)_{\lambda_0-n}= 0$  for $n\in \Z_{> 0}$. By Corollary \ref{dimen}, we have $\lambda_0=0$. Moreover, $(M/N)_{\lambda_0}$ has a 1-dimensional $A(L_{G(3)}(1, 0))$-submodule.

We now let $\C v$ be a 1-dimensional $A(L_{G(3)}(1, 0))$-submodule of $(M/N)_{\lambda_0}$, and $W$ be an $L_{G(3)}(1, 0)$-submodule of $M/N$ generated by $v$. We next prove that $W$ is an irreducible $L_{G(3)}(1, 0)$-module. Otherwise, let $U$ be a proper $L_{G(3)}(1, 0)$-submodule of $W$. Note that $U$ has the following decomposition $$U=\bigoplus_{\lambda\in\C}
U_{\lambda},$$ where
$U_\lambda=\{w\in U|(L(0)-\lambda)w=0\}$ and $\dim U_{\lambda}<\infty$ for any $\lambda\in \C$. Let $\mu_0$ be a number such that $U_{\mu_0}\neq 0$ but $U_{\mu_0-n}=0$  for $n\in \Z_{> 0}$. By Corollary \ref{dimen}, we have $\mu_0=0$. This forces that $U=W$, a contradiction. Thus, $W$ is an irreducible $L_{G(3)}(1, 0)$-module. In particular, by Proposition \ref{classify}, $W$ is isomorphic to $L_{G(3)}(1, 0)$.

We let $M^0$ be a submodule of $M$ such that $M^0$ contains $N$ and $M^0/N$ is isomorphic to $W$. By Lemma \ref{split2}, $M^0$ is isomorphic to $N\oplus L_{G(3)}(1, 0)$, this is a contradiction. Therefore, we have $M=N$ and $M$ is completely reducible.
\qed

\vskip.5cm
To summarize, by Theorems \ref{main}, \ref{semisimple}, $L_{G(3)}(1, 0)$ is a vertex operator superalgebra such that\\
 (i) $L_{G(3)}(1, 0)$ is not $C_2$-cofinite;\\
 (ii) The category of $L_{G(3)}(1, 0)$-modules is semisimple.

\end{document}